\documentclass{article}

\usepackage{amssymb, theorem, amsmath, rotate, multirow,epsfig}
\usepackage[all,2cell,dvips]{xy} \UseAllTwocells \SilentMatrices
\usepackage{fullpage}

\newcounter{nootje}
\newcommand\luijkabstract[1]
    {\begin{center}
     \small
     \begin{minipage}{4in}#1\end{minipage}
     \end{center}
    }
\newcommand\myref[1]{[#1]}
\newtheorem{theorem}{Theorem}[section]
\newtheorem{proposition}[theorem]{Proposition}

\newtheorem{lemma}[theorem]{Lemma}

\newtheorem{corollary}[theorem]{Corollary}
\newcommand\eindproof{\unskip\nobreak\hfill\hbox{\quad $\square$}\par \medskip}
\newenvironment{proof}{\noindent {\bf Proof.}}{\eindproof}
\newenvironment{proofof}{\noindent {\bf Proof of }}{\eindproof}
 
{\theorembodyfont{\rmfamily}

\newtheorem{remark}[theorem]{Remark}}

\makeatletter
\def\eqalign#1{\null\,\vcenter{\openup\jot\m@th
  \ialign{\strut\hfil$\displaystyle{##}$&$\displaystyle{{}##}$\hfil
      \crcr#1\crcr}}\,}
\makeatother

\newcommand\Lin{\mathop{\rm Lin} \nolimits}
\newcommand\Aut{\mathop{\rm Aut} \nolimits}
\newcommand\NS{\mathop{\rm NS} \nolimits}
\newcommand\Pic{\mathop{\rm Pic} \nolimits}
\newcommand\Spec{\mathop{\rm Spec} \nolimits}

\newcommand\rk{\mathop{\rm rank} \nolimits}
\newcommand\disc{\mathop{\rm disc} \nolimits}
\newcommand\Id{\mathop{\rm Id} \nolimits}
\newcommand\Q{\mathbb{Q}}
\newcommand\Qbar{\overline{\mathbb{Q}}}
\newcommand\R{\mathbb{R}}
\newcommand\Z{\mathbb{Z}}
\newcommand\F{\mathbb{F}}
\newcommand\p{\mathfrak{p}}
\renewcommand\P{\mathbb{P}}
\newcommand\Xbar{\overline{X}}
\renewcommand\O{\mathcal{O}}
\newcommand\kbar{\overline{k}}

\newcommand\Xint{\mathfrak{X}}
\newcommand\Ct{\tilde{C}}
\newcommand\isom{\cong}
\newcommand\et{\text{\rm \'et}}
\newcommand\ra{\rightarrow}

\newcommand\myrefart[6]{\item[{[#1]}] #2, #3, {\em #4}, #5, pp. #6.}
\newcommand\myrefartfive[5]{\item[{[#1]}] #2, #3, #4, #5.}
\newcommand\myrefartfour[4]{\item[{[#1]}] #2, #3, #4.}
\newcommand\myrefbook[5]{\item[{[#1]}] #2, {\em #3}, #4, #5.}

\newcommand\BT{BT}
\newcommand\chang{Ch}
\newcommand\friedman{Fr}
\newcommand\fulton{Fu}
\newcommand\hag{Ha}
\newcommand\kiran{Ki}
\newcommand\kedlaya{Ke}
\newcommand\matmon{MM}
\newcommand\ogus{NO}
\newcommand\pooneneuler{Po1}
\newcommand\poonen{Po2}
\newcommand\silvtwo{Si}
\newcommand\tate{Ta}
\newcommand\vangeemen{VG}

\newcommand\luijkheron{VL1}
\newcommand\luijkpic{VL2}

\begin{document}

\frenchspacing

\vspace*{-11mm}

\begin{center}
\sc {\bf \large Quartic K3 surfaces without nontrivial automorphisms}\\
 Ronald van Luijk
\end{center}

\vspace{-2mm}

\luijkabstract{For any field $k$ of characteristic at most $5$ we
  exhibit an explicit smooth quartic surface in $\P^3_k$ with trivial
  automorphism group over $\kbar$. We also show how this can be
  extended to higher  
  characteristics. Over $\Q$ we construct an example on which the set 
  of rational points is Zariski dense.}

\section{Introduction}


For any field $k$ we fix an algebraic closure of $k$, denoted by $\kbar$. 
For a variety $X\subset \P^{n+1}_k$ we set $\Xbar = X \times_k \kbar$, 
we say that $X$ is smooth if $\Xbar$ is regular, and we let $\Aut \Xbar$
denote the group of $\kbar$-automorphisms of $\Xbar$, while $\Lin \Xbar$
denotes the group of linear automorphisms of $\Xbar$, i.e., automorphisms
induced by a linear transformation of the coordinates of $\P^{n+1}$. 
The following theorem was proved by Poonen, see \myref{\poonen}, 
Thm. 1.6.

\begin{theorem}\label{poonen}
For any field $k$ and integers $n \geq 1$, $d \geq 3$ with $(n,d) \neq
(1,3)$, there exists a smooth hypersurface $X$ in $\P^{n+1}_k$ of degree $d$ 
with $\Lin \Xbar = \{1\}$.
\end{theorem}

\begin{remark}
For $(n,d) = (1,3)$ no such hypersurface exists. In that
case we get a plane cubic curve. With a flex as the origin,
such a curve obtains the structure of
an elliptic curve on which multiplication by $-1$ is a nontrivial
linear automorphism. 
\end{remark}

Poonen's proof of Theorem \ref{poonen} consists of giving carefully
crafted explicit examples. For fields that are large enough, a
nonconstructive proof can be deduced from various older results. For
details we refer to \myref{\poonen} and the references given there. 
The next theorem states that in many cases all automorphisms are
linear, see \myref{\chang}, Thm. 1, for dimension $1$, and
\myref{\matmon}, Thm. 2, for higher dimension.

\begin{theorem}\label{autislin}
If $X$ is a smooth hypersurface in $\P^{n+1}$ of degree $d$ with
$n\geq 1$, $d \geq 3$, and $(n,d) \not \in \{(1,3),(2,4)\}$, then
we have $\Aut \Xbar = \Lin \Xbar$.
\end{theorem}

Combining Theorem \ref{poonen} and \ref{autislin} we find the
following corollary.
\begin{corollary}\label{poonen2}
For any field $k$ and integers $n \geq 1$, $d \geq 3$ with $(n,d) \not \in
\{(1,3),(2,4)\}$, there exists a smooth hypersurface $X$ in
$\P^{n+1}_k$ of degree $d$ with $\Aut \Xbar = \{1\}$.
\end{corollary}

In this paper we deal with the remaining case $(n,d)=(2,4)$ in small
characteristics, including zero. The following theorem
states our main result.  
\begin{theorem}\label{maintheorem}
Let $k$ be any field of characteristic at most $5$. Then there exists a
smooth quartic surface $X \subset\P^3_k$ with $\Aut \Xbar = \{1\}$.
\end{theorem}

Smooth quartic surfaces in $\P^3$ are examples of K3 surfaces. Some
do have nonlinear automorphisms. The essential difference with the
hypersurfaces in Theorem \ref{autislin} is that the canonical sheaf is
trivial and the Picard group $\Pic \Xbar$ may be larger than $\Z$. 
The arithmetic of K3 surfaces is not well understood in general. It is for
instance not known whether there exists a K3 surface over a number
field on which the set of rational points is neither empty nor dense. 
Bogomolov and Tschinkel \myref{\BT} proved that if a K3 surface $X$ over a
number field $K$ has an infinite automorphism group or it admits an
elliptic fibration (not necessarily with a section), then there exists
a finite field extension $L$ of 
$K$, such that the set $X(L)$ of $L$-rational points on $X$ is Zariski
dense. The next theorem shows that it is not always necessary to
extend $K$, even if the automorphism group is trivial and we have
$K=\Q$.

\begin{theorem}\label{mainzariski}
Let $X \subset\P^3_{\Q}(x,y,z,w)$ be the surface given by 
$$
\eqalign{
2x^3w + 2x^2z^2 &+ 2x^2zw + x^2w^2 + 2xy^2w + 6xyz^2 + 2xyw^2 + xz^2w + \cr
   &+2y^3z + 6y^2z^2 + y^2w^2 + 10yz^2w + z^3w + 7z^2w^2 - 4zw^3=0.\cr
}
$$
Then $X$ is smooth, we have $\Aut \Xbar = \{1\}$, and the set
$X(\Q)$ is Zariski dense in $X$.  
\end{theorem}

The surface in Theorem \ref{mainzariski} admits an elliptic 
fibration. If Bogomolov and Tschinkel's result
does not apply, then there may still be infinitely many rational
points, as shown by the next theorem.

\begin{theorem}\label{maininfinite}
Let $X \subset\P^3_{\Q}(x,y,z,w)$ be the surface given by 
$$
w(3x^3+3xy^2-xyw+3xzw-xw^2+y^3-y^2w+2z^3+w^3) = (xy+xz+yz)(xy+xz+3yz).
$$
Then $X$ is smooth and does not admit an elliptic fibration over
$\Qbar$. We have $\Aut \Xbar = \{1\}$ and the set $X(\Q)$ is infinite.  
\end{theorem}

To prove Theorem \ref{maintheorem} we will write down explicit
examples of quartic surfaces $X$ and show that we have 
both $\Aut \Xbar = \Lin \Xbar$ and $\Lin \Xbar =\{1\}$. This is done in 
section \ref{withline} and in a different way in section
\ref{withconic}. The main idea 
behind both these sections is described in section \ref{secttheidea}.
As in \myref{\poonen}, it is not particularly hard to check that our
examples have no nontrivial linear automorphisms. 
The hard part was to find examples for 
which it is doable to verify this by hand. In section \ref{mainproof}
we prove Theorems \ref{maintheorem}, \ref{mainzariski}, and
\ref{maininfinite}. 

Our method works for any characteristic. We are limited, however, by
the lack of the ability to compute the characteristic polynomial
of Frobenius acting on certain cohomology groups in large characteristics. 
This is used to show $\Aut \Xbar = \Lin \Xbar$.

The author thanks Bjorn Poonen for bringing this problem to his
attention, Kiran Kedlaya for computing various characteristic
polynomials, and CRM at Montr\'eal for the hospitality and support
during the time in which this paper was written.

\section{The idea}\label{secttheidea}

The following lemma is the key ingredient
to proving $\Aut \Xbar = \Lin \Xbar$. 

\begin{lemma}\label{fromP3}
Let $X$ be a normal complete intersection in $\P^n$ 
and let $\sigma$ be a 
$\kbar$-automorphism of $\Xbar$. If $\sigma$ sends a hyperplane
section of $\Xbar$ to 
another hyperplane section of $\Xbar$, then we have $\sigma \in \Lin \Xbar$. 
\end{lemma}
\begin{proof}
Since $\sigma$ sends a hyperplane section to a hyperplane section, it
fixes $\O_X(1)$, so it sends a basis of $H^0(X,\O_X(1))$ to another
basis. By \myref{\hag}, Exc. II.8.4.c, the map $H^0(\P^n,\O_{\P^n}(1))
\rightarrow H^0(X,\O_X(1))$ is surjective, so this change of basis of
$H^0(X,\O_X(1))$ is induced by a change of basis of $H^0(\P^n,\O_{\P^n}(1))$.
The lemma now follows from \myref{\hag}, Thm. II.7.1.
\end{proof}

For the remaining of this paper we will assume that $X$ is a smooth
quartic surface in $\P^3$.
The condition of Lemma \ref{fromP3} is equivalent to $\sigma$ fixing
the class of hyperplane sections in the Picard group $\Pic \Xbar$.
This is the group of divisor classes on $\Xbar$ modulo linear
equivalence. The N\'eron-Severi group $\NS(\Xbar)$ is the group of
divisor classes modulo algebraic equivalence. For a precise definition of 
algebraic equivalence, see \myref{\hag}, Exc. V.1.7. The group
$\NS(\Xbar)$ is a finitely generated group. Its rank $\rho = \dim_\Q
\NS(\Xbar) \otimes \Q$ 
is called the Picard number of $X$. Note that in other papers 
$\rho$ is sometimes called the {\em geometric} Picard number of $X$.

For K3 surfaces, in particular for smooth quartic
surfaces in $\P^3$, linear equivalence is the same as algebraic and
numerical equivalence. This means that two divisors on $\Xbar$ are linearly
equivalent if and only if they have the same intersection number with
all other divisors. It implies that $\Pic \Xbar$ is finitely generated
and free, isomorphic to $\NS(\Xbar)$ and the intersection pairing endows 
$\Pic \Xbar$ with the structure of a lattice. For a divisor 
$D$, let $[D]$ denote the divisor class in $\Pic \Xbar \isom \NS(\Xbar)$. 
We will say that a divisor $D$ on $\Xbar$ is very ample if there exists an
integer $m$ and a closed immersion $\varphi \colon \Xbar \rightarrow
\P^m$ such that $\varphi(D)$ is  
a hyperplane section on $\varphi(\Xbar)$.

\begin{lemma}\label{adjunction}
Let $X$ be a smooth quartic surface in $\P^3$ and $C$ a smooth,
geometrically irreducible curve on $X$. Then we have $C^2 = 2g-2$,
where $g$ denotes the genus of $C$.
\end{lemma}
\begin{proof}
Since the canonical divisor on $X$ is trivial (see \myref{\hag},
Exm. II.8.20.3), this follows from the adjunction formula, see
\myref{\hag}, Prop. V.1.5.
\end{proof}

\begin{proposition}\label{theidea}
Let $X$ be a smooth quartic surface in $\P^3$ and let
$\sigma$ be a $\kbar$-automorphism of $\Xbar$. Let $H$ denote the divisor class
of hyperplane sections. Then the following implications hold. 
\begin{itemize}
\item[{\rm (1)}] If $\NS(\Xbar)$ is generated by $H$ then we have $\sigma
  \in \Lin \Xbar$.
\item[{\rm (2)}] If $\NS(\Xbar)$ is generated by $H$ and the divisor class
  of a line $L$, then we have $\sigma
  \in \Lin \Xbar$ and $\sigma$ fixes $L$. 
\item[{\rm (3)}] If $\NS(\Xbar)$ is generated by $H$ and the divisor class
  of a conic $C$, then we have $\sigma
  \in \Lin \Xbar$ and $\sigma$ fixes the plane that contains $C$.
\end{itemize}
\end{proposition}
\begin{proof}
Without loss of generality we may assume $\Xbar = X$.
Since the induced automorphism $\sigma_*$ of $\NS(\Xbar)$ preserves
intersection numbers, it is an automorphism of 
$\NS(\Xbar)$ as a lattice. Set $H'= \sigma_*(H)$.
Note that $\sigma_*$ sends very
ample divisor classes to very ample divisor classes and effective
divisor classes to effective divisor classes, so
$H'$ is a very ample divisor class. 

Case (1). The divisor classes $H'$ and $H$ both generate
$\NS(\Xbar)\isom \Z$, so we have either $H' = -H$ or $H' = H$. As $H'$
is effective, we find 
$H'=H$ and by Lemma \ref{fromP3} we have $\sigma \in \Lin \Xbar$.

Case (2). We have $H^2 = \deg X=4$ and $H\cdot [L] = \deg L = 1$, 
see \myref{\hag}, Exc. V.1.2. By Lemma \ref{adjunction} we have $L^2 = -2$.
This means that with respect to the basis $\{H,[L]\}$ the lattice $\NS(\Xbar)$ 
has Gram matrix 
$$
\left(
\begin{array}{cc}
4 & 1 \\
1 & -2 \\
\end{array}
\right)
$$
and thus discriminant $-9$.
The automorphism group of this lattice is isomorphic to $(\Z/2\Z)^2$
and it is generated by the automorphisms $[-1] \colon x \mapsto -x$ and 
$\tau \colon H \mapsto H+[L], [L] \mapsto -[L]$. The only automorphism
that sends effective divisor classes to effective divisor classes is
the identity, so we find $\sigma_* = \Id$. This implies that every
hyperplane section is sent to a divisor that is linearly equivalent to
it, i.e., another hyperplane section. From Lemma \ref{fromP3} 
we get $\sigma \in \Lin \Xbar$. We also conclude that $\sigma$ maps
$L$ to an effective divisor $L'$ that is linearly equivalent to $L$. 
As two different irreducible effective divisors have a nonnegative
intersection number, we conclude from $L' \cdot L = L^2 = -2$ that we
have $L'=L$, so $\sigma$ fixes $L$.

Case (3). We find similarly to case (2) that with respect to the basis
$\{H,[C]\}$ the lattice $\NS(\Xbar)$ has Gram matrix 
$$
\left(
\begin{array}{cc}
4 & 2 \\
2 & -2 \\
\end{array}
\right)
$$
and thus discriminant $-12$. 
Note that the conic $C$ is contained
in a plane $V$, see \myref{\hag}, Exc. IV.3.4. 
The other component $\Ct$ in $V\cap \Xbar$
has degree $\deg \Xbar \cdot \deg V - \deg C =2$, so it is also a
conic, which a priori might be degenerated.  
Let $a$ and $b$ be integers such that $H' = aH+b[C]$. 
Since $H'$ is very ample, it has positive intersection number with the
effective divisor classes $[C]$ and $[\Ct]=H-[C]$ (see \myref{\hag},
Thm. V.1.10), so we find 
$$
\eqalign{
0&< H' \cdot [C] = (aH+b[C])\cdot [C] = 2a-2b \cr
0&< H' \cdot [\Ct] = (aH+b[C])\cdot (H-[C]) = 2a+4b \cr
}
$$
We have $4=H^2=H'^2=4a^2+4ab-2b^2$, so we deduce 
$4 = (2a-2b)(2a+4b) + 6b^2 > 6b^2$, which implies $b=0$.
From $H'^2=4$ and the inequalities above we conclude $a=1$, so
$\sigma_* H = H'=H$.  
Again by Lemma \ref{fromP3} we have $\sigma \in \Lin \Xbar$. 
The orthocomplement of $H$ in $\NS(\Xbar)$ is generated by 
$D = 2[C]-H$, so we find $\sigma_* D = \pm D$. This implies 
$\sigma_*[C] = [C]$ or $\sigma_*[C] = H-[C] = [\Ct]$.
This means that $\sigma$ sends $C$ to a divisor that is linearly
equivalent to $C$ or to $\Ct$. As both $C$ and $\Ct$ have negative
self intersection, this implies as in case (2) that $\sigma$ maps $C$ to
$C$ or to $\Ct$. Since $V$ is the unique plane containing $C$ or $\Ct$,
the automorphism $\sigma$ fixes $V$. 
\end{proof}

The following lemma will be useful in conjunction with Proposition
\ref{theidea}. 

\begin{lemma}\label{genNS}
Let $X$ be a smooth quartic surface in $\P^3$ and let
$H$ denote the class of hyperplane sections. 
Let $\rho$ denote the Picard number of $X$.
Then the following implications hold. 
\begin{itemize}
\item[{\rm (1)}] If $\rho \leq 1$, then $\NS(\Xbar)$ is generated by $H$.
\item[{\rm (2)}] If $\rho \leq 2$ and $\Xbar$ contains a line $L$, then
  $\NS(\Xbar) = \langle H,[L] \rangle$.
\item[{\rm (3)}] If $\rho \leq 2$ and $\Xbar$ contains a conic $C$, then
  $\NS(\Xbar) = \langle H, [C] \rangle$.
\end{itemize}
\end{lemma}
\begin{proof}
Without loss of generality we may assume that the ground field is
algebraically closed, so we have $X = \Xbar$. We say that a lattice
$\Lambda$ with its pairing given by $(x,y)\mapsto x \cdot y$ is even
if the norm $x\cdot x$ is even for all $x \in \Lambda$. From Lemma
\ref{adjunction} we find that the integral lattice $\NS(X)$ is
generated by elements of even norm, so $\NS(X)$ is an even lattice. 
Note that the discriminants of a lattice $\Lambda$ and a sublattice 
$\Lambda'$ of finite index in $\Lambda$ are related by $\disc \Lambda' =
[\Lambda:\Lambda']^2  \cdot \disc \Lambda$.

Case (1). From $H^2=4$ we find $H\neq 0$, so $\rho=1$ and the lattice
$\langle H \rangle $ has finite index in $\NS(X)$. This implies 
$[\NS(X) : \langle H \rangle]^2 \cdot \disc \NS(X) = \disc \langle H
\rangle = H^2 = 4$. As any $1$-dimensional even lattice has even
discriminant, we find that $\disc \NS(X)$ is even, so $[\NS(X) :
  \langle H \rangle]=1$, and $\NS(X)$ is generated by $H$.

Case (2). As in the proof of Proposition \ref{theidea}, the lattice
$\langle H,[L]\rangle$ is $2$-dimensional and has discriminant
$-9$. Therefore we have $\rho=2$ and  
$[\NS(X) : \langle H, [L] \rangle]^2 \cdot \disc \NS = -9$. We conclude
$\NS(X) = \langle H, [L] \rangle$ or $\disc \NS(X) =-1$. By the
classification of even unimodular lattices, the latter case implies
that $\NS(X)$ is isomorphic to the lattice with Gram matrix 
$$
\left(
\begin{array}{cc}
0 & 1 \\
1 & 0 \\
\end{array}
\right),
$$
which is impossible for a quartic surface by a theorem of Van Geemen, see
\myref{\vangeemen}, 5.4. We find that $\NS(X)$ is indeed generated by
$H$ and $[L]$. 

Case (3). As in the proof of Proposition \ref{theidea}, the lattice
$\langle H,[C]\rangle$ is $2$-dimensional with discriminant $-12$ and 
Gram matrix 
\begin{equation}\label{twelve}
\left(
\begin{array}{cc}
4 & 2 \\
2 & -2 \\
\end{array}
\right). 
\end{equation}
Therefore we have  
$[\NS(X) : \langle H, [C] \rangle]^2 \cdot \disc \NS = -12$. We
conclude that the index $[\NS(X) : \langle H, [C] \rangle]$ divides
$2$. Take $D \in \NS(X)$. Then we have $2D = aH+b[C]$ for some
integers $a,b$. Since $\NS(X)$ is even, we find $8|4D^2 = (AH+b[C])^2
= 4a^2+4ab-2b^2$. This implies that $a$ and $b$ are both even, so we
have $D \in \langle H, [C] \rangle$. Since this holds for all $D \in
\NS(X)$, we find $\NS(X) = \langle H, [C] \rangle$. 
\end{proof}

In view of implication (1) of Proposition \ref{theidea} and Lemma
\ref{genNS}, one approach to proving Theorem 
\ref{maintheorem} is to take Poonen's examples of quartic surfaces $X$
with $\Lin \Xbar = \{1\}$ (see \myref{\poonen}) and prove that one of
them has Picard number $1$. There are at least two problems with this
approach. First 
of all, by Tate's conjecture (see \myref{\tate}) the N\'eron-Severi
group of a smooth quartic surface in $\P^3$ over a field of positive
characteristic has even rank. As Tate's conjecture has been proved for
all smooth quartic surfaces in characteristic $p\geq 5$ that are ordinary (see
\myref{\ogus}), there is not much hope for this approach in
positive characteristic. The second problem lies in proving that the 
Picard number of an explicit quartic surface $\Xbar$ over a
field of characteristic $0$ equals $1$.
The only way known to do this is described in \myref{\luijkpic}. It
requires two primes of good reduction for which we know the 
discriminant of the N\'eron-Severi lattice of the reduction
up to a square factor. It is not clear how to obtain this. 

Another way to use (1) of Proposition \ref{theidea} is to take the
quartic surfaces $X$ with Picard 
number $1$ that the author found in \myref{\luijkpic} and prove that
they satisfy 
$\Lin \Xbar = \{1\}$. The problem with this approach is that these
surfaces are defined by fairly erratic polynomials, which makes it hard
to get a handle on their linear automorphisms.

We will therefore use only implications (2) and (3) of Proposition
\ref{theidea}. This gives the extra advantage of knowing  
a specific subvariety that is fixed under any automorphism 
$\sigma$, which
%
%
severely restricts the possibilities for such $\sigma$. 

Let $X$ be a smooth surface over a number field $K$ and let $\p$ be a prime
of good reduction with residue field $k$.  
Let $\Xint$ be an integral model for $X$ over the localization $\O_\p$
of the ring of integers $\O$ of $K$ at $\p$. Let $l$ be any extension
field of $k$.  
Then by abuse of notation we will write 
$X \times l$ for $\Xint \times_{\Spec \O_\p} \Spec l$. 
To bound the Picard number of $X$ we will use the 
following propositions.  

\begin{proposition}\label{NSinj}
Let $X$ be a smooth surface over a number field $K$ and let $\p$ be a prime
of good reduction with residue field $k$. Then we have $$\rk \NS(X
\times \overline{K})
\leq \rk \NS(X \times \kbar).$$ 
\end{proposition}
\begin{proof}
See \myref{\fulton}, Example 20.3.6. It also follows from a natural injection 
$\NS(X \times \overline{K}) \otimes \Q \hookrightarrow \NS(X \times
\kbar) \otimes \Q$, see \myref{\luijkheron}, Prop. 6.2.
\end{proof}


\begin{proposition}\label{boundNS}
Let $X$ be a smooth surface over a finite field $k$ with $q$
elements. Let $l$ be a prime not dividing $q$. Let $F$ denote the 
automorphism on $H^2_{\et}(\Xbar,\Q_l)(1)$ induced by $q$-th power
Frobenius. Let $t$ denote the number of eigenvalues of $F$ that are
roots of unity, counted with multiplicity. Then there is a natural
injection 
$$
\NS(\Xbar) \otimes \Q_l \hookrightarrow H^2_{\et}(\Xbar,\Q_l)(1)
$$
that respects the action of Frobenius and we have $\rk \NS(\Xbar) \leq t$. 
\end{proposition}
\begin{proof}
See \myref{\luijkheron}, Prop. 6.2 and Cor. 6.4. Note that in the referred
corollary Frobenius acts on the 
cohomology group $H^2_{\et}(\Xbar,\Q_l)$ without a twist. Therefore,
the eigenvalues are scaled by a factor $q$. 
\end{proof}

\begin{remark}\label{tate}
Tate's conjecture (see \myref{\tate}) states that the inequality in
Proposition \ref{boundNS} is in fact an equality. 
\end{remark}

One way to find the characteristic polynomial of Frobenius is through
the Lefschetz formula. For a smooth, projective variety $X$ over the
field $\F_q$ and a positive integer $n$, 
this formula relates the number of points on $X$ over $\F_{q^n}$ to
the traces of $F_i^n$ for $0\leq i \leq 4$, where $F_i$ is the $q$-th
power Frobenius acting on 
$H^i_{\et}(\Xbar,\Q_l)$. For a K3 surface only the action for $i=2$
is unknown, so the trace of $F_2^n$ can be computed from
$\#X(\F_{q^n})$. From these traces for various $n$ one can deduce the
characteristic polynomial of $F_2$. 
If a priori a Galois invariant subspace $V$ of $H^2_{\et}(\Xbar,\Q_l)$
is known, then it suffices to find the characteristic polynomial of
Frobenius on the quotient, which requires fewer traces. In our case
$V$ will always be generated by the class of hyperplanes and the class
of either a line or a conic. For details see \myref{\luijkheron}, section
6 and 7, and \myref{\luijkpic}, section 2. 
A faster way to find the number $t$ as in
Proposition \ref{boundNS} is to use Kedlaya's algorithm based on De
Rham cohomology \myref{\kedlaya}, analogous to his algorithm for
hyperelliptic curves 
described in \myref{\kiran}. We will use both methods freely without
reproducing the details that can be found in these references. 

The following lemmas will be used in sections \ref{withline} and
\ref{withconic}.

\begin{lemma}\label{reduced}
Let $Z$ be a smooth, irreducible hypersurface in $\P^n$ over an
algebraically closed field, with $n \geq 3$. Let $H$ be a
hyperplane not equal to $Z$ and let $C$ be an irreducible component of
the scheme theoretic intersection 
$H \cap Z$. Then $C$ is reduced and the intersection multiplicity of
$H$ and $Z$ along $C$ is equal to $1$.
\end{lemma}
\begin{proof}
Let $S = k[x_0,\ldots, x_n]$ denote the homogeneous coordinate ring of
$\P^n$. Suppose $Z$ is given by $f = 0$ for some $f \in S$. 
Without loss of generality we may assume that $H$ is given by
$x_n=0$. Then $H$ is isomorphic to $\P^{n-1}$ with homogeneous
coordinate ring $T = k[x_0,\ldots, x_{n-1}]$ and $H\cap Z$ is given by
the image $g$ of $f$ under the homomorphism $S \rightarrow T$ that
sends $x_n$ to $0$. Since $T$ is a unique factorization domain, the
component $C$ corresponds to an irreducible factor $h$ of $g$. 
We find from the definition of intersection multiplicity
(see \myref{\hag}, p. 53) that the
intersection multiplicity of $H$ and $Z$ along $C$ is equal to $1$ if
and only if the exponent of $h$ in $g$ equals $1$, which is the case
if and only if the component $C$ is reduced. Suppose the exponent of
$h$ in $g$ were at least $2$. Then we could
write $f = x_n q +h^2 p$ for some $q \in S$ and $p \in R$. This
implies that for every point on $C$, which is given by $x_n = h =0$,
we have $f=0$ and $\partial f /\partial x_i=0$ for $i \in
\{0,\ldots,n-1\}$. From $n \geq 3$ we conclude that $C$ has dimension
at least $1$. From the Projective Dimension Theorem (see \myref{\hag},
Thm. I.7.2) it follows that there is a point $P$ on $C$ that also 
satisfies  $\partial f/\partial x_n = 0$, which implies that $Z$ is
singular at $P$. The lemma follows from this contradiction.
\end{proof}

\begin{lemma}\label{allnodescusps}
Let $k$ be a field of characteristic different from $2$ and $3$. 
Let $X$ be a K3 surface over $k$ with Picard number at most $2$. Let $\pi
\colon \, X \ra \P^1$ be an elliptic fibration. Then all singular fibers are
irreducible curves with either a node or a cusp. Let $n$ and $c$
denote the number of nodal and cuspidal fibers respectively. Then we
have $n+2c = 24$.
\end{lemma}
\begin{proof}
Let $F$ be a fiber of $\pi$ with $r$ irreducible components. By
\myref{\silvtwo}, Prop. III.8.2, these components generate 
a sublattice of $\NS(\Xbar)$ of rank $r$ in which every element $z$
satisfies $z^2 \leq 0$. 
Any ample divisor has positive self intersection and therefore none of
its multiples is 
contained in this sublattice, so the Picard number of $\Xbar$ is at
least $r+1$. This proves that all fibers are irreducible. In general,
all irreducible singular fibers of an elliptic fibration
are nodal or cuspidal curves.
For any fiber $F$, let $e(F)$ denote the Euler characteristic of
$F$. For K3 surfaces, we have $\sum_F e(F)=24$, where the sum is taken
over all singular fibers $F$ of $\pi$, see \myref{\friedman},
Cor. 7.16 and p. 178. 
The lemma follows from the fact that outside characteristic $2$ and
$3$, nodal and cuspidal fibers 
have Euler characteristic $1$ and $2$ respectively, see
\myref{\pooneneuler}, p. 266.
\end{proof}

\begin{lemma}\label{func}
Let $f\in \R[x]$ be a polynomial of degree $d$ all of whose roots have complex
absolute value equal to $1$. Then we have $x^df(x^{-1}) = (-1)^nf(x)$,
where $n$ is the order of vanishing of $f$ at $x=1$.
\end{lemma}
\begin{proof}
The only real numbers of absolute value $1$ are $\pm 1$, so the
factorization of $f$ into irreducible factors over $\R$ is 
$$
f = c(x-1)^n(x+1)^m \prod_{i=1}^r g_i,
$$
for some constant $c$, integers $m,n,r$, and quadratic monic polynomials
$g_i$. Since the conjugate roots of $g_i$ have absolute value $1$, we
have $g_i = x^2+a_ix+1$ for some $a_i\in \R$, so $x^2g_i(x^{-1}) = g_i(x)$. 
The lemma follows. 
\end{proof}

\section{Explicit examples containing a line}\label{withline}

\begin{proposition}\label{exampleline}
Let $k$ be any field and let $f_1,f_2,f_3 \in k[x,y,z,w]$ be
homogeneous polynomials of degree $2$. Let $X$ be the surface in
$\P^3$ given by  
\begin{equation}\label{eqexline}
w(x^3+xy^2+wf_1)+z(y^3+zf_2)+zwf_3=0.
\end{equation}
Let $\pi \colon \, \Xbar \ra \P^1$ be the morphism given by
$[x:y:z:w]\mapsto [z:w]$. 
Suppose that $X$ is smooth and its Picard number is at most $2$.
Assume that the fiber of $\pi$ above $[1:0]$ is singular at
$[0:0:1:0]$, the fiber above $[0:1]$ has a cusp at
$[0:0:0:1]$, and no other fiber is cuspidal. Then we have
$\Aut \Xbar = \{1\}$. 
\end{proposition}

Let $L$ be the line given by $w=z=0$. 
The morphism $\pi$ in Proposition \ref{exampleline} is an
elliptic fibration (not necessarily with a section). Each fiber is
the union of the components other than $L$ in some hyperplane section 
through $L$. In order for $f_1$, $f_2$, and $f_3$ to satisfy all
conditions of Proposition \ref{exampleline}, they have to satisfy other
conditions that are easier to check. The following lemma states some
of these conditions. Not only are they useful for finding explicit
examples, some of them will also be used in the proof of
Proposition \ref{exampleline}. 

\begin{lemma}\label{wzcubed}
Let $f_1$, $f_2$, $f_3$, $X$, and $\pi$ be as in Proposition
\ref{exampleline}. Then there are $f_1',f_2',f_3'\in k[x,y,z,w]$ 
that satisfy the following conditions.
\begin{itemize}
\item[{\rm (a)}] $X$ is given by $w(x^3+xy^2+wf_1')+z(y^3+zf_2')+zwf_3'=0$,
\item[{\rm (b)}] $f_1',f_2' \in k[x,y]$,
\item[{\rm (c)}] $f_1'$ is a square over $\kbar$,
\item[{\rm (d)}] the coefficients of $z^2$ and $w^2$ in $f_3'$ are
  nonzero,
\item[{\rm (e)}] the polynomials $x^3+xy^2+wf_1'$ and $y^3+zf_2'$ are
  irreducible over $\kbar$, 
\item[{\rm (f)}] if $k$ is finite, say $\# k=q$, and $F$ is the $q$-th
  power Frobenius acting on $H^2_{\et}(\Xbar,\Q_l)(1)$, and Tate's
  conjecture is true, then the sign of the functional equation for the
  characteristic polynomial of $F$ is positive.
\end{itemize}
\end{lemma}
\begin{proof}
By definition there are $f_1',f_2',f_3'$ such that (a) is satisfied,
namely given by $f_i'=f_i$. 
After collecting all monomials of the polynomial 
$g=w(x^3+xy^2+wf_1')+z(y^3+zf_2')+zwf_3'$ that are divisible
by $zw$ in the term $zwf_3'$ we may assume $f_1' \in k[x,y,w]$ and $f_2'
\in k[x,y,z]$. 
Since the fiber $F_\infty$ of $\pi$ above $[1:0]$, given by
$w=y^3+zf_2' =0$, contains the
singular point $P_z = [0:0:1:0]$, we find that the polynomial 
$y^3+zf_2'$ and its derivatives with respect to $x$, $y$, and $z$ all
vanish at $P_z$. This implies that $f_2'$ does not contain the variable
$z$, so we have $f_2' \in k[x,y]$. Similarly
we find $f_1' \in k[x,y]$ from the singularity $P_w = [0:0:0:1]$ in
the fiber $F_0$ above $[0:1]$, which proves (b). The fact that $F_0$ 
has a cusp at $P_w$ is then equivalent to 
the fact that $f_1'$ is a square in $\kbar[x,y]$, which proves (c). 
It is now easily checked that $g$ and its derivatives
with respect to $x,y$, and $w$ all
vanish at $P_w$. Since $X$ is smooth at $P_w$, this implies that
$\partial g/\partial z$ does not vanish at $P_w$, which is equivalent
to the coefficient of $w^2$ in $f_3'$ being nonzero. Similarly we find
from the fact that $X$ is smooth at $P_z$ that the coefficient of
$z^2$ in $f_3'$ is nonzero. This takes care of (d). For (e),
note that all fibers are geometrically integral by Lemmas
\ref{reduced} and \ref{allnodescusps}. The fibers $F_0$ and $F_\infty$
are given by $z=x^3+xy^2+wf_1'=0$ and $w=y^3+zf_2'=0$
respectively, so the polynomials $y^3+zf_2'$ and
$x^3+xy^2+wf_1'$ are irreducible over $\kbar$. Finally, assume all
hypotheses of (f). By Lemma \ref{genNS} the N\'eron-Severi group $\NS(\Xbar)$
is generated by the class of hyperplane sections and the divisor class
of the line $L$. Since $F$ acts trivially on these
classes, by Proposition \ref{boundNS} and Remark \ref{tate} we find
that there are exactly 
two eigenvalues of $F$ that are roots of unity, both equal to $1$. 
In particular, the multiplicity of the eigenvalue $1$ is even,
which implies that the sign of the functional equation is positive by
Lemma \ref{func}.
\end{proof}

%

\begin{proofof}{\bf Proposition \ref{exampleline}.}
Replace the $f_i$ by the $f_i'$ of Lemma \ref{wzcubed}.
Since $X$ contains the line $L$, we find from Lemma \ref{genNS} that 
$\NS(\Xbar)$ is generated by the class of hyperplane sections and the
class $[L]$. Let $\sigma$ be any $\kbar$-automorphism of $\Xbar$.
By Proposition \ref{theidea} we have $\sigma \in
\Lin \Xbar$ and $\sigma$ fixes $L$. The first claim implies that there
exists a 
matrix $A$ such that $\sigma$ sends $[x:y:z:w]$ to $[x':y':z':w']$ 
with $(x'\,\,y'\,\,z'\,\,w')^t = A \cdot (x\,\,y\,\,z\,\,w)^t$, where
$v^t$ denotes the transpose of the vector $v$.
Since $A$ maps the generators $z,w$ of the ideal of $L$ to generators
of the same ideal, we find that $A$ is of the form 
$$
\left(
\begin{array}{cccc}
a & b & \kappa & \lambda \\ 
c & d & \mu & \nu \\
0 & 0 & p & q \\
0 & 0 & r & s \\
\end{array}
\right).
$$ 
Let $g$ be the polynomial in the left-hand side of equation
(\ref{eqexline}). Set $(x'\,\,y'\,\,z'\,\,w')^t = A \cdot
(x\,\,y\,\,z\,\,w)^t$. As $\sigma$ is an automorphism, 
the polynomial $g' = g(x',y',z',w')$ in terms of the
variables $x,y,z,w$ also defines $X$, so $g'$ is a 
scalar multiple of $g$. 
After scaling the matrix $A$ we may assume $g'=g$. 
Comparing the coefficients in $g$ and $g'$ of the monomials that
are linear in $z$ and $w$ we find that the following expressions are
all zero.
$$
\eqalign{
Q_1&=c^3p + (a^3 + ac^2)r, \cr
Q_2&=c^3q + (a^3 + ac^2)s - 1,\cr
Q_3&=3c^2dp + (3a^2b + 2acd + bc^2)r,\cr
Q_4&=3c^2dq + (3a^2b + 2acd + bc^2)s,\cr
Q_5&=3cd^2p + (3ab^2 + ad^2 + 2bcd)r,\cr
Q_6&=3cd^2q + (3ab^2 + ad^2 + 2bcd)s - 1,\cr
Q_7&=d^3p + (b^3 + bd^2)r - 1, \cr
Q_8&=d^3q + (b^3 + bd^2)s. \cr
}
$$
This implies that we also have 
$$
\eqalign{
cr &= -sdQ_3+rdQ_4+scQ_5-rcQ_6 =0, \cr
c^3s+d^3r &= sd^3Q_1-rd^3Q_2-sc^3Q_7+rc^3Q_8 = 0.\cr
}
$$
From the first of these last two equations we find $r=0$ or $c=0$. In case
of $r=0$, we find $s\neq 0$ from the fact that $A$ is invertible and
from the second equation we then find $c=0$, so we have $c=0$ in either
case. Since $A$ is invertible, we conclude $a\neq 0$. Then the
equation $Q_1=0$ implies $r=0$ and $Q_2=0$ gives $s=a^{-3}$. From
$Q_4=0$ it now follows that we have $3b=0$. Equation $Q_6=0$ yields
$s=a^{-1}d^{-2}$, which together with $s=a^{-3}$ gives $d = \pm a$.  
The equation $Q_7=0$ shows $p=d^{-3}$.

From $r=0$ we deduce that the fiber $F_\infty$ of $\pi$ above $[1:0]$ is
mapped to itself. All fibers are geometrically integral
by Lemma \ref{reduced} and \ref{allnodescusps}. As integral cubic plane
curves have at most one singular point (see \myref{\hag}, Exm. V.3.9.2),
the singular point $[0:0:1:0]$ of $F_\infty$ is unique and thus fixed by
$\sigma$. This implies $\kappa=\mu=0$. 
The cuspidal fiber $F_0$ above $[0:1]$ is sent isomorphically to the fiber 
above $[q:s]$. The only cuspidal fibers lie above $[0:1]$ and possibly
$[1:0]$, so from $s \neq 0$ we find $q=0$. We deduce that $\sigma$
sends the uniue cusp $[0:0:0:1]$ of $F_0$ to itself, which 
implies $\lambda =\nu = q =0$. 

From Lemma \ref{wzcubed} we know that 
$y^3+zf_2$ is irreducible, so $f_2 \in k[x,y]$ is not divisible by
$y$. Therefore the coefficient $\chi$ of $x^2$ in $f_2$ is nonzero.
Comparing the coefficient of $x^2z^2$ in $g'$ and $g$ we find from
$a^2=d^2$ that we have 
$$
\chi = \chi a^2p^2 = \chi a^2d^{-6}   = \chi a^{-4},
$$
which implies that we have $a^4 = 1$. By Lemma \ref{wzcubed} the
coefficient of 
$wz^3$ in $g$ is also nonzero. Comparing this to the coefficient of $wz^3$
in $g'$ we now find $1 = d^{-9}a^{-3} = d^{-1}a$, so $a=d$. 
After dividing $A$ by
$a$ we may assume $1=a=d=p=s$.  By comparing the coefficient of 
$xyz^2$ in $g$ and $g'$ we then get $2\chi b = 0$, so
$2b=0$. Together with $3b=0$ this gives $b=0$ and we find
that $\sigma$ is the identity. As this holds for every $\sigma \in
\Aut \Xbar$ we find $\Aut \Xbar = \{1\}$. 
\end{proofof}

\begin{remark}\label{howpicked}
The main idea of the proof of Proposition \ref{exampleline} is that
the line $L$ given by $z=w=0$ is fixed by any 
automorphism. This implies that any automorphism 
permutes the fibers of $\pi$. We picked some
extra conditions to ensure that some of these fibers have to be
fixed, which yields enough information to deduce that every
automorphism is trivial. These conditions are by no means general.  
They are carefully chosen such that on one hand they are easy enough
to verify the proof without a computer, while on the other hand it
leaves ample examples among which to search for surfaces with the
right Picard number.
\end{remark}

Before we give an explicit example that satisfies all conditions of
Proposition \ref{exampleline}, we will sketch the odds we had to beat
to find one. Based on Lemma \ref{wzcubed}, the least number of
monomials of an 
equation as in Proposition \ref{exampleline} that satisfies all
conditions of the proposition is $7$. 
This minimum is achieved if and only if
there are $\alpha,\beta,\gamma,\delta \in k^*$, such that the quartic
surface $X$ is given by 
\begin{equation}\label{minimal}
w(x^3+xy^2+\alpha wy^2) + z(y^3+\beta zx^2) +zw(\gamma z^2 +\delta
w^2) =0. 
\end{equation}
The map $[x:y:z:w] \mapsto [-x:-y:z:w]$ gives an isomorphism between
the surfaces corresponding to $(\alpha, \beta, \gamma, \delta)$ and 
$(-\alpha, -\beta, -\gamma, -\delta)$, while the map 
$[x:y:z:w] \mapsto [x:-y:-z:w]$ gives an isomorphism between
the surfaces corresponding to $(\alpha, \beta, \gamma, \delta)$ and 
$(\alpha, \beta, -\gamma, -\delta)$.
For $k=\F_3$ this gives four quadruples of
isomorphic surfaces. The four nonisomorphic surfaces are all smooth.
For one of them the fibration $\pi$ has reducible components, so its
Picard number is at least $3$ by 
Lemma \ref{allnodescusps}. For the remaining three surfaces $X$ we computed
the characteristic polynomial of Frobenius on
$H^2_{\et}(\Xbar,\Q_l)(1)$. The sign of the functional
equation of one of them is negative, which gives a bound for the
Picard number larger than $2$ by Lemma \ref{wzcubed}. Unfortunately,
Proposition \ref{boundNS}  
gives an upper bound of at least $4$ for the remaining two Picard
numbers as well. 

We therefore look at surfaces given by (\ref{minimal}) with an extra 
term $\epsilon yz^2w$ for some $\epsilon \in \F_3^*$. In this case we get 
eight nonisomorphic smooth surfaces, one of which has a hyperplane
section that contains three components. For only one of the remaining
seven surfaces the sign of the corresponding functional equation is
positive. We computed the characteristic polynomial of Frobenius in
that case and found from Proposition \ref{boundNS} that the Picard
number is at most $2$. This yields the following corollary.

\begin{corollary}\label{charthree}
Let $k$ be any field of characteristic $0$ or $3$. Let $X$ in $\P^3$
over $k$ be given by 
$$
x^3w - x^2z^2 + xy^2w + y^3z + y^2w^2 + yz^2w - z^3w + zw^3=0
%
%
$$
Then we have $\Aut \Xbar = \{1\}$.
\end{corollary}
\begin{proof}
Let $g$ denote the polynomial in the given equation.
Note that the equation is of the form used in Proposition
\ref{exampleline} with $f_1 = y^2$, $f_2=-x^2$, and $f_3 = yz-z^2+w^2$. 
Although without a computer algebra package it is a bit of work
to show that $X$ is smooth, with such a package one easily checks
that this is the case. 
Let $X_{0}$ and $X_3$ denote the surfaces defined by $g=0$ over $\Q$
and $\F_3$ respectively.
We will show $\rk \NS(\Xbar_3) \leq 2$. 
Each of the nonsingular fibers of the fibration $\pi$ as in
Proposition \ref{exampleline} can be given the structure of an
elliptic curve by searching for a rational point on it. 
There are fast algorithms to compute the number of points on
elliptic curves, implemented in for instance {\sc magma}. Summing over
all fibers, we find the number of points on $X_3$ over $\F_{3^n}$ 
for $n \in \{1,\ldots, 10\}$. These numbers are 
$$
17, 
101 ,
812 ,
6545 ,
58502 ,
531902 ,
4788164 ,
43074713 ,
387494444 ,
3486985076.
$$
As in \myref{\luijkheron}, Prop. 7.1, and \myref{\luijkpic},
Thm. 3.1, we can use the Lefschetz formula to find the characteristic
polynomial of Frobenius acting on $H^2_{\et}(\Xbar_3,
\Q_l)(1)$ from these numbers. Here we use the fact that we already
know a Galois invariant subspace of $H^2_{\et}(\Xbar_3,\Q_l)(1)$,
generated by a hyperplane 
section and the line $L$. In $\Q[t]$ the characteristic polynomial
factors into irreducible factors as 
$$
\eqalign{
\frac{1}{3}(t-1)^2
(3t^{20}&- t^{19}- t^{17}+ 2t^{16}+ 2t^{15}- t^{13}-t^{12}+ t^{11}\cr 
& + t^{10}+ t^{9}- t^{8}- t^{7}+ 2t^5+ 2t^{4}- t^{3}- t + 3).\cr
}
$$ 
The roots of the irreducible factor of degree $20$ are not integral,
so they are not roots of unity. By Proposition \ref{boundNS} we find 
$\rk \NS(\Xbar_3) \leq 2$. From Proposition \ref{NSinj} we also 
find $\rk \NS(\Xbar_{0}) \leq 2$. Depending on the characteristic we
have either $X = X_3 \times k$ or $X = X_{0}\times k$. As the
N\'eron-Severi group is algebraic, we conclude $\NS(\Xbar) =
\NS(\Xbar_3)$ or $\NS(\Xbar) = \NS(\Xbar_{0})$. In both cases we
conclude $\rk \NS(\Xbar) \leq 2$. The fibers of $\pi$ above 
$[0:1]$ and $[1:0]$ have a cusp at the points $[0:0:0:1]$ and $[0:0:1:0]$
respectively. In characteristic $0$ there are $20$ more singular
fibers. By Lemma \ref{allnodescusps} these fibers are all nodal curves. 
In characteristic $3$ there are $14$ more singular fibers.
A tedious calculation shows that again none of the corresponding
fibers has a cusp. 
From Proposition \ref{exampleline} we deduce $\Aut \Xbar = \{1\}$.
\end{proof}

\section{Explicit examples containing a conic}\label{withconic}

\begin{proposition}\label{examplemine}
Let $k$ be any field with elements $\alpha,\beta$ satisfying 
$\alpha^3\beta \neq \alpha\beta^3$. Let $f\in k[x,y,z,w]$ be a
homogeneous polynomial of degree $3$, such that the coefficients of
$y^3$ and $z^3$ in $f$ are different, or the coefficients of $x^2y$
and $x^2z$ in $f$ are different. Suppose that the surface $X$ in $\P^3$ 
given by 
$$
w f = (xy+xz+\alpha yz)(xy+xz+\beta yz)
$$
is smooth with Picard number at most $2$. Then we have $\Aut
\Xbar = \{1\}$. 
\end{proposition}

As in the previous section, the following lemma will be useful both
for the proof of this 
proposition and for constructing examples that satisfy all conditions.

\begin{lemma}\label{xyzcubed}
Let $f$ be as in Proposition \ref{examplemine}. Then the following
conditions hold. 
\begin{itemize}
\item[{\rm (a)}] The coefficients of $x^3$, $y^3$, and $z^3$ in $f$
  are nonzero.
\item[{\rm (b)}] If $k$ is finite, say $\# k=q$, and $F$ is the $q$-th
  power Frobenius acting on $H^2_{\et}(\Xbar,\Q_l)(1)$, 
  and Tate's conjecture is true, then the sign
  of the functional equation for the characteristic polynomial of $F$
  is positive.
\item[{\rm (c)}] Suppose $k'$ is a finite quadratic subfield of $k$, say
  $\# k'=q$, the elements $\alpha$ and $\beta$ are conjugate over $k'$, and $f
  \in k'[x,y,z,w]$. Then $X$ is defined over $k'$. If $F$ is the $q$-th  
  power Frobenius acting on $H^2_{\et}(\Xbar,\Q_l)(1)$, 
  and Tate's conjecture is true, then the sign of the functional
  equation for the characteristic polynomial of $F$ is negative.
\end{itemize}
\end{lemma}
\begin{proof}
Note that $X$ contains the points $P_1=[1:0:0:0]$, $P_2=[0:1:0:0]$,
and $P_3 = [0:0:1:0]$. One easily checks that for $i=1,2,3$, the
surface $X$ is smooth at $P_i$ if and only if we have $f(P_i) \neq
0$. This implies (a). For (b) and (c), note that by 
Lemma \ref{genNS} the N\'eron-Severi group $\NS(\Xbar)$
is generated by the class $H$ of hyperplane sections and the divisor class
of the conic $C$ given by $w=xy+xz+\alpha yz=0$. In case (b), $F$ acts
trivially on these classes, and the proof proceeds as the proof of (f)
of Lemma \ref{wzcubed}. In case (c), $F$ acts trivially on $H$, and
by a nontrivial quadratic character on the class $2[C]-H$. By
Proposition \ref{boundNS} and Remark \ref{tate}
this implies that counted with multiplicity, there are exactly two
eigenvalues of $F$ that are roots of unity, equal to $1$ and $-1$
respectively. In 
particular, the multiplicity of the eigenvalue $1$ is odd, which
implies that the sign of the functional equation is negative by Lemma
\ref{func}.
\end{proof}

\begin{proofof}{\bf Proposition \ref{examplemine}.}
Let $C$ and $\Ct$ denote the conics given by $w=xy+ xz+\alpha yz=0$ and
$w=xy+xz+\beta yz=0$ respectively. As both are contained
in $X$, we find from Lemma \ref{genNS} that $\NS(\Xbar)$ is generated
by the class of hyperplane sections and the class $[C]$. Let $\sigma$
be a $\kbar$-automorphism of $\Xbar$. From Proposition
\ref{theidea} we conclude that we have $\sigma \in \Lin \Xbar$ and
that $\sigma$ fixes the plane containing $C$, which also contains
$\Ct$. It follows that $\sigma$ permutes the intersection points of
$C$ and $\Ct$. These are $P_x=[1:0:0:0]$, $P_y=[0:1:0:0]$, and
$P_z=[0:0:1:0]$, where the first has multiplicity $2$ and the others have
multiplicity $1$. This implies that $\sigma$ fixes $P_x$ and either it
also fixes $P_y$ and $P_z$, or it switches these two. 
Therefore $\sigma$ is given by $[x:y:z:w] \mapsto [x':y':z':w']$ with 
$(x' \,\, y'\,\, z'\,\,w')^t = A (x \,\, y\,\, z\,\,w)^t$ for 
a matrix $A$ of the form 
$$
\left(
\begin{array}{cccc}
p & 0 & 0 & \lambda \\ 
0 & q & 0 & \mu \\
0 & 0 & r & \nu \\
0 & 0 & 0 & s \\
\end{array}
\right)
, \qquad \text{or} \qquad
\left(
\begin{array}{cccc}
p & 0 & 0 & \lambda \\ 
0 & 0 & r & \mu \\
0 & q & 0 & \nu \\
0 & 0 & 0 & s \\
\end{array}
\right),
$$
with $p,q,r,s \in \kbar^*$ and $\lambda, \mu, \nu \in \kbar$. Let $g$
denote the polynomial $w f - (xy+xz+\alpha yz)(xy+xz+\beta yz)$.

Suppose we are in the first (diagonal) case. Then 
with $x' = px+\lambda w$, $y' = qy+\mu w$, $z'=rz+\nu w$, and $w' =
sw$ the polynomial $g' = g(x',y',z',w')$ in terms of the variables
$x$, $y$, $z$, and $w$ also defines $\Xbar$, 
so it is a constant multiple of the polynomial $g$.
After scaling $A$ we may assume $g'=g$. 
For any monomial $M$, let $c_M$ and $c'_M$ denote the coefficients of
$M$ in $g$ and $g'$ respectively. Then for every monomial $M$ we have
$c_M = c_M'$. From the assumption $\alpha^3\beta \neq \alpha\beta^3$
we find $\alpha+\beta\neq 0$, so
$c_M \neq 0$ for $M \in S_1 = \{xy^2z, xyz^2, x^2y^2\}$. 
As the coefficient of $x^3$ in $f$ is nonzero by Lemma
\ref{xyzcubed}, we also conclude $c_{wx^3} \neq 0$. 
Note that the terms $(xy+xz+\alpha yz)(xy+xz+\beta yz)$ and $(x'y'+x'z'+\alpha
y'z')(x'y'+x'z'+\beta y'z')$ do not contribute to $c_M$ and $c'_M$ for 
$M = wx^3$. Therefore the equations
$c'_M=c_M$ for $M\in S_1 \cup \{wx^3\}$ give
$sp^3=pq^2r=pqr^2=p^2q^2=1$, which implies $p=q=r=s$. After scaling
$A$ we may assume $p=q=r=s=1$.
Then for $M \in \{xy^2w,xz^2w,yz^2w\}$ the equations $c_M= c_M'$ give 
$$
\eqalign{
-2\lambda -(\alpha+\beta)\nu +c_{xy^2w}&=c_{xy^2w}, \cr
-2\lambda -(\alpha+\beta)\mu +c_{xz^2w}&=c_{xz^2w}, \cr
-(\alpha+\beta)\lambda - 2\alpha\beta\mu+c_{yz^2w}& = c_{yz^2w}. \cr
}
$$
By hypothesis we have $\alpha^2 \neq \beta^2$. One easily checks that 
this implies that this system of linear equations has no nontrivial
solutions in $\lambda, \mu$, and $\nu$, so we have $\lambda=\mu=\nu=0$
and $\sigma$ is the identity.

Suppose we are in the second (nondiagonal) case. Now we have  $x' =
px+\lambda w$, $y' = rz+\mu w$, $z'=qy+\nu w$, and $w' = sw$. As
before we may assume that $g' = g(x',y',z',w')$ equals $g$, and we
can again compare the coefficients $c_M$ and $c_M'$ of the monomials
$M$ in $g$ and $g'$ respectively.
As in the first case we find $p=q=r=s$, so we can scale $A$ to get
$p=q=r=s=1$. This gives $c_{wz^3}'=c_{wy^3}$. From
$c_{wz^3}'=c_{wz^3}$ we then deduce $c_{wy^3}=c_{wz^3}$, so the
coefficients of $y^3$ and $z^3$ in $f$ are equal.
For $M \in \{x^2yw,x^2zw\}$ the equations $c_M=c_M'$ give 
$$
\eqalign{
-2\mu-2\nu + c_{x^2zw} & = c_{x^2yw}, \cr
-2\mu-2\nu + c_{x^2yw} & = c_{x^2zw}, \cr
}
$$
which has no solution in any characteristic unless $c_{x^2yw} =
c_{x^2zw}$. This implies that also the coefficients of $x^2z$ and
$x^2y$ in $f$ are equal, which contradicts our assumptions. We
conclude that this case does not occur, so 
$\sigma$ is the identity. As this holds for every $\sigma \in
\Aut \Xbar$ we find $\Aut \Xbar = \{1\}$. 
\end{proofof}

\begin{remark}
As for Proposition \ref{exampleline}, the conditions in Proposition
\ref{examplemine} are by no means general.  
\end{remark}

\begin{remark}
We will give examples of quartic surfaces $X$ for which
Proposition \ref{examplemine} immediately tells us that we have 
$\Aut \Xbar = \{1\}$. To do this over a field $k$  
we need  $\alpha, \beta \in \kbar$, such
that $\alpha\beta(\beta-\alpha)(\beta+\alpha)\neq 0$ and the
right-hand side of the equation in Proposition \ref{examplemine} is
defined over $k$. If $k$ has at least $4$ elements than the existence
of such $\alpha$ and $\beta$ follows 
from the fact that for any nonzero $\alpha$ the polynomial $\alpha
X(X^2-\alpha^2)$ has at most $3$ roots. For $k=\F_2$ we can take 
$(\alpha,\beta) = (\zeta,\zeta+1)$ for some $\zeta \in \F_4$ with 
$\zeta^2+\zeta+1=0$. For $k=\F_3$ we can take $(\alpha,\beta) =
(1+i,1-i)$ for some element $i \in \F_9$ with $i^2 = -1$.
\end{remark}


\begin{corollary}\label{chartwo}
Let $k$ be any field of characteristic $0$ or $2$ and let $X\subset
\P^3$ over $k$ be given by
$$
w(x^3+y^3+z^3+x^2z+xw^2)=x^2y^2+2x^2yz+x^2z^2-xy^2z-xyz^2+y^2z^2.
$$
Then we have $\Aut \Xbar = \{1\}$.
\end{corollary}
\begin{proof}
Let $X_0$ and $X_2$ be the surfaces defined by the given equation over
$\Q$ and $\F_2$ respectively. 
Although it is tedious work without a
computer algebra package, one easily checks that $X_0$ and $X_2$ are
smooth. That means that $p=2$ is a prime of good reduction for $X_0$. 
In \myref{\luijkpic}, Rem. 6, it was shown that we have $\rk
\NS(\Xbar_2) \leq 2$. 
From Proposition \ref{NSinj} we also find $\rk \NS(\Xbar_0) \leq 2$.
Because the N\'eron-Severi group is algebraic, we find as in the proof
of Corollary \ref{charthree} that we have $\rk \NS(\Xbar) \leq 2$. 
Let $\zeta$ be an element (possibly 
in a quadratic extension of $k$) that satisfies $\zeta^2+\zeta+1=0$. 
Then the equation for $X$ is of the form given in Proposition
\ref{examplemine} with 
$$
(\alpha,\beta,f)=(\zeta, \zeta^2, x^3+y^3+z^3+x^2z+xw^2).
$$
From Proposition \ref{examplemine} we find $\Aut \Xbar = \{1\}$.
\end{proof}

\section{Proof of the main theorems}\label{mainproof}

For characteristics $0$, $2$, and $3$ we have seen examples of quartic
surfaces with trivial automorphism group in Corollary \ref{charthree}
and \ref{chartwo}. Kiran Kedlaya has also found an example of a triple
$(\alpha, \beta,f)$ in characteristic $5$ that satisfies all
conditions of Proposition \ref{examplemine}.

\begin{proposition}\label{charfive}
Take $k= \F_5$, $\alpha = 1$, $\beta=3$, and 
$$
f = 3x^3+3xy^2-xyw+3xzw-xw^2+y^3-y^2w+2z^3+w^3.
$$
Then all conditions of Proposition \ref{examplemine} are satisfied.
\end{proposition}
\begin{proof}
The bound of the Picard number is the only condition that is
not trivial to check, see \myref{\kedlaya}. 
\end{proof}

\begin{proofof}{\bf Theorem \ref{maintheorem}.}
Corollary \ref{chartwo} gives an explicit example 
in characteristic $0$ and $2$. For characteristic $3$ an example is 
given by Corollary \ref{charthree}. For characteristic $5$ we use 
Proposition \ref{examplemine} with $(\alpha,\beta,f)$ as
in Proposition \ref{charfive}.
\end{proofof}

\begin{proofof}{\bf Theorem \ref{mainzariski}.}
Let $f$ denote the polynomial in the given equation and set
$f_1 = \frac{1}{2}(x+y)^2$, $f_2 = x^2+3xy+3y^2$, and 
$f_3 = x^2+\frac{1}{2}xz+5yz+\frac{1}{2}z^2+\frac{7}{2}zw-2w^2$. Then we have 
$\frac{1}{2}f = w(x^3+xy^2+wf_1)+z(y^3+zf_2)+zwf_3$, so $X$ is as in
Proposition \ref{exampleline}. One easily checks that $X$ is smooth
with good reduction at $3$. Let $X_3$ denote the reduction at $3$. 
For $n = 1,\ldots, 10,$ we compute the number of points on
$X_3$ over $\F_{3^n}$. These numbers are 
$$
15,107 ,639 ,6935 ,59790 ,533729 ,4790661 ,43039079 ,387592263, 3486831422.
$$
As in the proof of Corollary \ref{charthree}, and as in
\myref{\luijkheron}, Prop. 7.1, and \myref{\luijkpic},
Thm. 3.1, we can find the characteristic
polynomial of Frobenius acting on $H^2_{\et}(\Xbar_3,
\Q_l)(1)$ from these numbers. In $\Q[x]$ the characteristic polynomial
factors into irreducible factors as
$$
\eqalign{
\frac{1}{3}(x-1)^2(
3x^{20} &+ x^{19} - x^{18} + 5x^{17} - 3x^{15} + 4x^{14} - 2x^{13} -
        2x^{12} + x^{11} \cr
 &- 3x^{10} + x^9 - 2x^8 - 2x^7 + 4x^6 - 3x^5 + 5x^3 - x^2 + x + 3).\cr
}
$$
As the roots of the factor of degree $20$ are not integral,
we find $\rk \NS(\Xbar_3) \leq 2$ from Proposition \ref{boundNS}. 
Then from Proposition \ref{NSinj} we also find $\rk \NS(\Xbar) \leq 2$.
As before, let $\pi$ denote the elliptic fibration given by 
$[x:y:z:w]\mapsto [z:w]$.
The fiber $F_{\infty}$ above $[1:0]$ has a node at $[0:0:1:0]$. 
The fiber $F_0$ above $[0:1]$ has a cusp at $[0:0:0:1]$. 
There are $21$ more singular fibers, all conjugate over $\Q$. 
By Lemma \ref{allnodescusps} these are all nodal curves.
We conclude that all conditions of Proposition \ref{exampleline} are
fulfilled, so we find $\Aut \Xbar = \{1\}$.

Consider the curve $C = X \cap H$, where $H$ is the hyperplane given
by $w+x=z$. Note that $C$ can be parametrized by $\tau \colon\,\P^1 \ra C,
\, [s:1] \mapsto [z-w:y:z:w]$ with 
$$
y= - 8s^4 - 4s^3 +6s^2 -1, \qquad z=2s^2(4s^2+2s-1),
\qquad w= 2s.
$$
The curve $C$ intersects each fiber $F$ of $\pi$ with
multiplicity $\deg F =3$, so the restriction of $\pi$ to $C$ has
degree $3$. Identify the function field $k(C)$ of $C$ with $k(s)$
through $\tau$. Identify the function field of the base curve $\P^1$
of $\pi$ with $k(t)$. The composition $\pi\tau$ is given by 
$[s:1] \mapsto [s(4s^2+2s-1):1]$, so the corresponding field extension  
$k(s)/k(t)$ of degree $3$ is given by $t = s(4s^2+2s-1)$. The map 
$\tau$ corresponds to a point $\O$, defined over $k(s)$, on the
generic fiber $E$ of $\pi$, giving $E$ the structure of an elliptic
curve over $k(s)$. The sum $S$ of the two conjugates of $\O$, both 
defined over $\overline{k(t)}$, is also defined over $k(s)$. 
The map $\pi\tau \colon \, \P^1(s) \ra \P^1(t)$ is ramified at
$s_0=\frac{1}{6}$ above $t_0 = -\frac{5}{54}$, where by abuse of
notation we denote the point $[\alpha : \beta]$ by $\alpha\beta^{-1}$.
The third point above
$t_0$ is $s_1 = -\frac{5}{6}$. Specializing at $s = s_1$ we find 
that $E$ specializes to the elliptic fiber $F_{t_0}$ above $t_0$, the
point $\O$ specializes to $\O_{s_1} = \tau(s_1) = [295: 263: 25:
  -270]$, and the point $S$ specializes to $S_{s_1} = 2Q$ for $Q =
\tau(s_0) = [59: 139: 5: -54]$. 
With these explicit equations, a standard computation shows that 
$Q$, and therefore $S_{s_1}$, has infinite order on the elliptic
curve $F_{t_0}$ with $\O_{s_1}$ as origin. This implies that $S$ has 
infinite order on the generic fiber $E$. The infinitely many multiples of
$S$ correspond to infinitely many rational curves on $X$, all with
infinitely many rational points, obtained from specializing $s$. 
This shows that the set $X(\Q)$ of rational points is Zariski dense in $X$.
\end{proofof}

\begin{remark}
In the notation of \myref{\BT}, the curve $C$ in the proof of Theorem
\ref{mainzariski} is a rational {\em nt}-multisection of $\pi$. 
This follows from the proof. A highbrow argument for this fact is that
$C$ is saliently ramified, because the restriction of $\pi$ to $C$
ramifies in the smooth fiber above $t=-\frac{5}{54}$, see \myref{\BT},
Prop. 4.4. 

To find a surface with such a multisection, we fixed a hyperplane $H$, 
three points $P_1$, $P_2$, and $P_3$, and searched for 
$f_1$, $f_2$, and $f_3$ such that for $X$ as in Proposition
\ref{exampleline} the intersection $C = X \cap H$ is singular at the
points $P_1$, $P_2$, and $P_3$. As the arithmetic genus of any
hyperplane section is $\frac{1}{2}(4-1)(4-2)=3$, this ensures that $C$ has
geometric genus $0$, see \myref{\hag}, Exm. V.3.9.2. Since $C$
intersects the line $L$ given by $w=z=0$ in a smooth point, we can
parametrize $C$. 
Per point $P_i$ we get three equations for the
indeterminate coefficients of the $f_j$. If the $P_i$ are collinear,
then $X \cap H$ will be the union of a line and a cubic curve, which
means that the Picard number of $X$ will be more than $2$. We
therefore pick the $P_i$ such that modulo $3$ they are not collinear. 
Setting $f_1 = a(y+bx)^2$, $f_2 = cx^2+dxy +ey^2$ and $f_3 = \sum_M
c_M M$, where $M$ runs over all $10$ monomials of degree $2$ in the
variables $x,y,z,w$, we find $9$ equations in $15$ variables. Choosing 
$H$ to be given by $w+x=z$, and 
$[1:1:0:-1]$, $[-1:-1:1:2]$, and $[1:-1:1:0]$ for the $P_i$, the
equations reduce to $9$ linearly independent linear equations, so
modulo $3$ we find $3^{15-9} = 729$ candidate surfaces. Out of all
these surfaces, only 157 are nonsingular. For 64 of these, the morphism
$\pi$ has a reducible fiber, making the Picard number more than
$2$. Out of the remaining 93 surfaces, 70 have a reducible hyperplane
section defined over $\F_9$ that does not contain $L$, which also
implies larger Picard number than $2$. For the 
last $23$ surfaces we computed the first $10$ traces of powers of
Frobenius. For ten of these surfaces the sign of the functional
equation could be determined to be negative, so the bound for the
Picard number is larger than $2$ by Lemma \ref{wzcubed}.
For eight surfaces the sign could not be determined yet. For four
of these, either sign would give an upper bound of at least
$4$. Determining the sign for the other four requires substantial 
extra computing time.
For the last five surfaces we could determine the sign to be positive. 
For only one of these, the
upper bound for the Picard number is $2$. This is the surface of
Theorem \ref{mainzariski}.
\end{remark}

\begin{proofof}{\bf Theorem \ref{maininfinite}.}
The given surface reduces modulo $5$ to the surface of Proposition
\ref{charfive}. From that Proposition we find that $X$ has good reduction 
at $5$ and since the reduction has Picard number $2$, we also find
$\rk \NS(\Xbar) \leq 2$ from Proposition \ref{NSinj}. All conditions
of Proposition \ref{examplemine} are fulfilled, so we find 
$\Aut \Xbar =\{1\}$. By Proposition \ref{genNS} the
N\'eron-Severi group is generated by the class $H$ of hyperplane sections
and the class of the conic $C$ given by $w=xy+xz+yz=0$, with $H^2=4$,
$C^2=-2$, and $H\cdot C = 2$. For any nonzero class $D = aH+b[C]$ we find 
$D^2 = (2a+b)^2-3b^2$, so $D^2\neq 0$. Because a fiber of any
fibration of $X$ has self intersection $0$, we conclude that there is
no such fibration. The conic $C$ has a point $[1:0:0:0]$, so it
contains infinitely many rational points. 
\end{proofof}

\section*{References}

\begin{itemize}
\myrefart{\BT}{Bogomolov, F. and Tschinkel, Yu.}{Density of rational
  points on elliptic K3 surfaces}{Asian J. Math.}{{\bf 4}, 2
  (2000)}{351--368}
\myrefart{\chang}{Chang, H.C.}{On plane algebraic curves}{Chinese
  J. Math.}{{\bf 6 (2)} (1978)}{185--189}
\myrefbook{\friedman}{Friedman, R.}{Algebraic surfaces and holomorphic
vector bundles}{Universitext, Springer}{1998}
\myrefbook{\fulton}{Fulton, W.}{Intersection  Theory, second
  edition}{Springer}{1998} 
\myrefbook{\hag}{Hartshorne, R.}{Algebraic Geometry}{GTM {\bf 52},
Springer-Verlag, New-York}{1977}
\myrefart{\kiran}{Kedlaya, K.}{Counting points on hyperelliptic curves
    using Monsky-Washnitzer cohomology}{J. Ramanujan Math. Soc.}{{\bf
    16} (2001),  no. 4}{323--338}
\myrefartfour{\kedlaya}{Kedlaya, K.}{Bounding Picard numbers of
  surfaces using $p$-adic cohomology}{{\em in preparation}}
\myrefart{\matmon}{Matsumura, H. and Monsky, P.}{On the automorphisms
  of hypersurfaces}{J. Math. Kyoto Univ.}{{\bf 3} (1963/1964)}{347--361}
\myrefart{\ogus}{Nygaard, N. and Ogus, A.}{Tate's conjecture for K3
  surfaces of finite height}{Ann. of Math.}{{\bf 122}
  (1985)}{461--507}
\myrefart{\pooneneuler}{Poonen, B.}{An explicit algebraic family of
  genus-one curves violating the Hasse principle}{J. de Th\'eorie
  des Nombres de Bordeaux}{{\bf 13} (2001)}{263--274}  
\myrefart{\poonen}{Poonen, B.}{Varieties without extra automorphisms
  III: hypersurfaces}{Finite Fields Appl.}{{\bf 11} (2005)}{230--268}
\myrefbook{\silvtwo}{Silverman, J.H.}{Advanced Topics in the Arithmetic of
Elliptic Curves}{GTM {\bf 151}, Springer-Verlag, New-York}{1994}
\myrefart{\tate}{Tate, J.}{Algebraic cycles and poles of zeta
  functions}{Arithmetical Algebraic Geometry}
  {ed. O.F.G. Schilling (1965)}{93--110}
\myrefart{\vangeemen}{Van Geemen, B.}{Some remarks on Brauer groups of K3
  surfaces}{Advances in Math.}{{\bf 197} (2005)}{222--247}
\myrefartfive{\luijkheron}{Van Luijk, R.}{An elliptic K3 surface
associated to Heron triangles}{{\em preprint}, \\ 
  available at {\tt arXiv:math.AG/0411606}}{2004} 
\myrefartfive{\luijkpic}{Van Luijk, R.}{K3 surfaces with Picard number
  one and infinitely many rational points}{{\em preprint}, available at
 {\tt arXiv:math.AG/0506416}}{2005} 
\end{itemize}

\end{document}